\documentclass[]{article} 
 
\usepackage{amsmath}[1996/11/01]
\usepackage{amssymb,amsfonts}

\DeclareMathOperator{\cols}{cols}

\DeclareMathOperator{\Sym}{Sym}

\DeclareMathOperator{\Hom}{Hom}

\DeclareMathOperator{\End}{End}
\DeclareMathOperator{\GL}{GL}
\DeclareMathOperator{\Aut}{Aut}
\DeclareMathOperator{\SL}{SL}
\DeclareMathOperator{\St}{St}

\DeclareMathOperator{\Cl}{Cl}

\DeclareMathOperator{\N}{N}
\DeclareMathOperator{\nr}{nr}

\advance \textheight by -1 \topmargin
\advance \textheight by -1 \headheight
\advance \textheight by -1 \headsep
\advance \textheight by -1 \footskip
\vsize = \textheight
\hsize = \textwidth
\def\CC{{\Bbb C}}

\def\HH{{\Bbb H}}

\def\OO{{\Bbb O}}

\def\QQ{{\Bbb Q}}
\def\RR{{\Bbb R}}

\def\ZZ{{\Bbb Z}}

\newtheorem{theorem}{Theorem}[section]

\newcommand{\bew}{\noindent\underline{Proof.}\ }
\newtheorem{remark}[theorem]{Remark}

\newtheorem{lemma}[theorem]{Lemma}
\newtheorem{proposition}[theorem]{Proposition}
\newtheorem{definition}[theorem]{Definition}

\newtheorem{folg}[theorem]{Corollary}
\newtheorem{corollary}[theorem]{Corollary}

\renewcommand{\OO}{{\mathcal{O}}}
\newcommand{\trace}{\mbox{trace}}

\newcommand{\eb}{\phantom{zzz}\hfill{$\square $}\smallskip}
\renewcommand{\a}{\mathfrak{a}}

\title{Maximal finite subgroups 
and minimal classes}
\author{Renaud Coulangeon and Gabriele Nebe}

\begin{document}

\maketitle

\begin{abstract}
We apply Voronoi's algorithm to compute 
representatives of the conjugacy classes of maximal finite subgroups 
of the unit group of a maximal order in some simple $\QQ $-algebra.
This may be used to 
show in small cases that non-conjugate orders have non-isomorphic 
unit groups.
\end{abstract}

\section{Introduction}
Let $A$ be some simple $\QQ $-algebra and let 
$\Lambda $ and $\Gamma $ be two maximal orders in $A$. 
If $A$ is not a division algebra, then the 
order $\Lambda $ is generated by its unit group $\Lambda ^{\times}$ as 
a $\ZZ $-lattice (see Lemma \ref{genunits}). So $\Lambda ^{\times}$ and $\Gamma ^{\times}$ are conjugate in 
$A^{\times}$ if and only if the two orders $\Lambda $ and $\Gamma $ are conjugate, 
which can be decided with the arithmetic theory of orders exposed in the next
section.
 By the 
theorem of Skolem and Noether we hence have that 
 the unit groups are 
\emph{conjugate}  if and only if $\Lambda $ and $\Gamma $ 
are isomorphic as orders over the center of $A$.
The motivation of this paper is to develop tools for deciding 
whether the two unit groups are \emph{isomorphic}, which is in general much more difficult than the conjugacy problem. 
In fact this innocent question was raised by Oliver Braun 
during his work on the paper \cite{Braun} that 
grew out of his Bachelor thesis in Aachen supervised by the second author. 

One invariant of the isomorphism class of 
$\Lambda ^{\times}$ is  
the number of conjugacy classes of maximal finite 
subgroups. 
Our main result is that these maximal finite subgroups arise 
as automorphism groups of \emph{well rounded minimal classes}, 
which will be defined in Section \ref{minclass}. 
The basic idea underlying this approach is already apparent in Ry{\v{s}}kov's paper 
\cite{Rysh72} on the computation of the  finite subgroups of $\GL_n(\ZZ)$.
 Nevertheless, whereas Ry{\v{s}}kov's classify all finite subgroups and then develop
 \emph{ad hoc} arguments to determine the maximal ones,
 our method permits in principle to solve the problem directly. 
Precisely, a refinement of the classical Voronoi algorithm, 
involving Berg\'e-Martinet-Sigrist's equivariant version
 of Voronoi's theory \cite{BMS92},
is applied to compute the cellular decomposition of 
a suitable retract of a cone of positive definite Hermitian forms, and therewith also 
the (finitely many) conjugacy classes
of maximal finite subgroups of $\Lambda ^{\times} $.  As will be illustrated in Section \ref{exples},
 this turns out to be enough, in some cases,
 to distinguish between non-isomorphic unit groups.
 The argument can of course not be reversed: non-isomorphic unit groups might 
have the same conjugacy classes of maximal finite subgroups.
 Note also that, as in the classical case of $\GL_n(\ZZ)$, 
the obtained cellular decomposition 
can be used to
compute the integral homology of $\Lambda ^{\times}$. The relevance of Voronoi theory in such homology computations was first highlighted in the works of Soul\'e \cite{Sou75,Sou78} and Ash \cite{Ash80,Ash84}, and it has given rise since then to numerous developments (we refer the interested reader to P. Gunnels' appendix of \cite{Ste07} which provides an excellent survey on this topic, and to \cite{DES11,Rah13,RF11} for recent related works, especially on Bianchi groups).

The methods apply to arbitrary (semi)-simple $\QQ $-algebras, 
though we are mainly interested in the case 
 where 
$A$ is a matrix ring over either a imaginary quadratic 
number field or a definite rational quaternion algebra. 
For these algebras we may ease these computations by adopting a
projective notion of minimal vectors as exposed  in
Section \ref{mini}.

\section{Conjugacy classes of maximal orders}\label{conjclass} 
The theory in this section is well known and can be extracted from the two books 
\cite{Reiner} and \cite{Deuring}. 
However, we did not find a self-contained short exposition of the proof 
of Theorem \ref{conj}, so we repeat the details here for the reader's convenience. 
Let $A$ be a simple $\QQ $-algebra.
Then $A = M_n(K)$ for 
some rational division algebra $K$ with center $Z(K)$.
Let $R$ be the maximal order in $Z(K)$ and choose 
some maximal $R$-order $\OO $ in $K$.
An $\OO $-lattice $L$ of rank $n$ is a finitely generated $\OO $-submodule
of the right $K$-module $V:=K^n$ that contains a $K$-basis. 
By  Steinitz-theorem
(see for instance \cite[Theorem 4.13, Corollary 35.11]{Reiner}) there 
are right ideals $\mathfrak{c}_1, \ldots 
\mathfrak{c}_n$ of $\OO $ 
 and a basis
$(e_1,\ldots , e_n)$ of $V$
such that
$$L=e_1 \mathfrak{c}_1  \oplus ... \oplus e_n \mathfrak{c}_n .$$
The family $\left( \mathfrak{c}_i, e_i\right)_{1\leq i \leq n}$ is called a \emph{pseudo-basis} of $L$. 
The \emph{Steinitz-invariant} of $L$, denoted $\St(L)$, is the class 
$$\St(L):= \left[ \mathfrak{c}_1\right]  +\cdots +\left[ \mathfrak{c} _n\right] $$ in the 
group $\Cl(\OO )$ of stable isomorphism classes of right $\OO $-ideals and does not depend of the choice of a pseudo-basis.
By Eichler's theorem (see \cite[Theorem (35.14)]{Reiner}) the reduced norm 
$$\nr : \Cl(\OO ) \to \Cl _K(R) $$ 
induces a group isomorphism between 
$\Cl(\OO )$ and the ray class group $\Cl_K(R)$, the quotient of the ideal group of $R$
modulo those principal ideals $\alpha R$ for which $\iota (\alpha ) > 0$ for all
real places $\iota $ of $Z(K)$ that ramify in $K$. 

If $n\geq 2$ (which we assume in the following) then, 
as a consequence of Corollary 35.13 of \cite{Reiner},
two lattices $L_1,L_2 \leq V$ are isomorphic as $\OO $-modules, if and only if 
they have the same Steinitz-invariant. In particular,
$L $ is isomorphic to $L({\mathfrak c})$ where 
$$L({\mathfrak c}) = e_1 \OO \oplus \ldots \oplus e_{n-1} \OO  \oplus e_n \mathfrak{c} $$
for any ideal $\mathfrak{c} $ with $[\mathfrak{c} ] = \St(L)$.
The endomorphism ring
$$\End _{\OO }(L) = \left\lbrace  X\in M_n(K) \mid X L \subseteq L \right\rbrace  $$
is a maximal order in $\End _K (V) \cong A$. In fact any maximal
order in $A$ is obtained this way (see \cite[Corollary 27.6]{Reiner}).
If $[\mathfrak{c} ] = \St(L)$ then $\End_{\OO }(L)$ is conjugate in
$\GL_n(K)$ to
$$\End _{\OO }(L(\mathfrak{c})) = \Lambda (\mathfrak{c}  ) := \left( \begin{array}{cccc}
\OO  & \ldots & \OO  & \mathfrak{c}^{-1} \\ 
\vdots & \ldots & \vdots & \vdots \\
\OO  & \ldots & \OO  & \mathfrak{c}^{-1} \\ 
\mathfrak{c} & \ldots & \mathfrak{c} & \OO '  \end{array} \right) $$
where $\OO ' = O_{l}(\mathfrak{c} ) = \left\lbrace  x\in K \mid x \mathfrak{c} \subseteq \mathfrak{c} \right\rbrace  $.

\begin{lemma} \label{genunits} 
For $n\geq 2$ any maximal order $\Lambda $ in $A=M_n(K)$ is generated as
a $\ZZ $-order by its unit group.
\end{lemma}

\bew
Without loss of generality let $\Lambda = \Lambda (\mathfrak{c}) $ and let
$$(x_1,\ldots , x_d) , \ (y_1,\ldots , y_d), \ (z_1,\ldots , z_d) $$
be $\ZZ $-bases of $\OO $, $\mathfrak{c} $, respectively $\mathfrak{c}^{-1}$. 
We denote by $e_{ij}$ the matrix units in $M_n(K)$ having an entry $1$ at $i,j $ and 
$0$ elsewhere, and $I_n = e_{11}+\ldots + e_{nn} $ the unit matrix. 
Let $X $ be the $\ZZ $-order spanned by $\Lambda (\mathfrak{c})^{\times}$.
Since $I_n$ and $I_n+x_k e_{ij} $ $\in \Lambda (\mathfrak{c})^{\times}$ 
we obtain that $x_k e_{ij} \in X$
 for all 
$k=1,\ldots , d$, $1\leq i\neq j \leq n-1 $.
Similarly $y_k e_{ni}$ and $z_je_{in} $, as well as $y_kz_j e_{nn} $ and 
$z_jy_ke_{ii} $ are in $X$ for all $i=1,\ldots n-1$, $k,j = 1,\ldots , d$. 
As the $y_kz_j$ generate $\OO '$ and the $z_jy_k$ generate $\OO $ the order 
$X$ contains $\Lambda (\mathfrak{c})$.
\eb

\begin{corollary} \label{conunits}
Let $\Lambda $ and $\Gamma $ be two maximal orders in the simple algebra $A$ and
assume that $A$ is not a division algebra.
Then $\Lambda ^{\times}$ and $\Gamma ^{\times}$ are conjugate in $A^{\times}$ if and only if 
$\Lambda $ and $\Gamma $ are conjugate. 
\end{corollary} 

A separating invariant of the conjugacy classes of maximal orders in $A$ 
can be constructed in a suitable class group of the center of $A$. 

\begin{definition}
Let $\Cl_K (n) := \Cl_K(R) / \langle  \nr({\mathfrak a}) ^n \mid \mathfrak{a} \unlhd \OO \rangle $ denote the quotient of the ray class group $\Cl _K(R)$ 
defined above modulo the $n$-th powers of the reduced norms of the two-sided $\OO $-ideals. 
\end{definition} 

Note that the subgroup $\langle  \nr({\mathfrak a}) ^n \mid \mathfrak{a} \unlhd \OO \rangle $
can be obtained from the discriminant of $K$. In particular it does not depend on the 
choice of the maximal order $\OO $. 
Also if $K$ is commutative then $\Cl_K(n) = \Cl(K) / \Cl(K)^n$ is just the 
class group of $K$ modulo the $n$-th powers. 

\begin{theorem} \label{conj} Let $A= M_n(K)$ be a simple $\QQ$- algebra and $\OO$ a maximal order in $K$. For any two right $\OO$-ideals $\mathfrak{c}$ and  $\mathfrak{c}'$, the corresponding maximal orders
$\Lambda (\mathfrak{c}) $ and $\Lambda (\mathfrak{c}')$ are conjugate in $A^{\times}=\GL_n(K)$ 
if and only if $\nr ([\mathfrak{c} ]) = \nr([\mathfrak{c}'] )$ in $\Cl_K (n) $.
\end{theorem} 

\bew
We use the approach in \cite[Section VI.8]{Deuring}. 
Let $\Gamma := M_n(\OO ) = \Lambda (\OO )$. 
Then any other maximal order in $A$ arises as the left order of some $\Gamma $-right ideal,
in particular 
$$\Lambda (\mathfrak{c} ) = O_l (I({\mathfrak c})) = \left\lbrace  a\in A \mid aI({\mathfrak c}) \subseteq 
I({\mathfrak c}) \right\rbrace  
\mbox{ where } I({\mathfrak c}) = \left( \begin{array}{ccc} \OO & \ldots & \OO \\ 
\vdots & \ldots & \vdots \\ 
\OO & \ldots & \OO \\ 
\mathfrak{c} & \ldots & \mathfrak{c} \end{array} \right) .$$ 
Two left orders $O_l(I)$ and $O_l(I')$ are conjugate, if and only if 
$I' = a I J$ for some $a\in A^{\times}$ and some two-sided 
fractional $\Gamma $-ideal $J$.  By Morita theory 
any two-sided  $\Gamma $-ideal $J$ is of the form
$J = M_n({\mathfrak a}) $ for some two-sided $\OO $-ideal ${\mathfrak a} $ in $K$. 
By \cite[Lemma (35.8)]{Reiner}, the reduced norm of $J=\Hom_{\OO } (\OO ^n , {\mathfrak a}^n)$
 equals 
$\nr ({\mathfrak a})^n \in \Cl_K(R) $ and the reduced norm of $I({\mathfrak c}) = 
\Hom _{\OO } (L(\OO ),L({\mathfrak c}) ) $ is $\nr ({\mathfrak c})$.
By \cite[Theorem 35.14]{Reiner} the reduced norm 
is injective, 
so 
$$ 
I({\mathfrak c}) = a I({\mathfrak c} ') M_n({\mathfrak a}) \mbox{ for some } a\in A^{\times} \ 
\mbox{ if and only if } \ \nr ({\mathfrak c}) = \nr ({\mathfrak c'}) \nr ({\mathfrak a})^n . 
$$
\eb

\section{Positive cones} 
Let $K$ be some rational division algebra and $A=M_n(K)$. 
Then $A_{\RR }:= A \otimes _{\QQ } \RR $ 
is a semi-simple real algebra, 
hence a direct sum of matrix rings over one of $\HH $, $\RR $ or $\CC $.
It carries a canonical involution that we use to define symmetric elements.
Let $d$ denote the degree of $K$, so $d^2 = \dim _{Z(K)} (K)$, and let
$$ \begin{array}{ll} 
\iota _1 ,\ldots , \iota _s  &   \mbox{ be 
the real places of $Z(K)$ that ramify in $K$,} \\
 \sigma _1,\ldots , \sigma _r &
\mbox{ the real places of $Z(K)$ that do not ramify in $K$ } \\
\tau _1,\ldots , \tau _t  & 
\mbox{ the complex places of $Z(K)$.} 
\end{array}
$$
Then
$$K_{\RR }:= K\otimes _{\QQ } \RR \cong 
\bigoplus _{i=1}^{s} M_{d/2}(\HH ) \oplus 
\bigoplus _{i=1}^{r} M_d(\RR ) \oplus 
\bigoplus _{i=1}^{t} M_d( \CC )  .$$
The ``canonical'' involution  $\ ^*$ (depending on the choice of this
isomorphism)
is defined on 
any simple summand of $K_{\RR }$
to be transposition for $M_d(\RR )$, transposition and complex 
(respectively quaternionic) conjugation for $M_d(\CC )$ and 
$M_{d/2} (\HH )$. 
The resulting involution  on $K _{\RR }$ is again denoted by $\ ^*$.
As usual it defines 
a mapping $\ ^{\dagger} : M_{m,n} (K_{\RR } ) \to M_{n,m}(K_{\RR } )$
by applying $\ ^*$ to the entries and then 
transposing the $m\times n$-matrices. 
In particular this defines an involution $\ ^{\dagger} $ on 
$A _{\RR} = M_n(K_{\RR}) $.
In general this  involution  will not 
fix the set $A$.

\begin{definition}
$\Sigma:= \Sym (A _{\RR }) :=  \left\lbrace  F\in A_{\RR} \mid F^{\dagger } = F \right\rbrace  $ 
is the $\RR $-linear subspace of symmetric elements of $A _{\RR} $.
It supports the positive definite inner product 
$$\langle F_1,F_2 \rangle := \trace (F_1F_2 ) $$ 
where $\trace $ is the reduced trace of the semi simple $\RR $-algebra $A _{\RR}$. The  real vector space $\Sigma $  contains the open real cone of positive elements 
$$ {\mathcal P} := \left\lbrace   (q_1,\ldots , q_s, f_1,\ldots , f_r, h_1,\ldots , h_t) \in \Sigma \mid q_i, f_j ,h_k \mbox{ pos. def.} 
 \right\rbrace  .$$
\end{definition}

Let $V$ be the simple left $A $-module $K^{n}$. 
Then $V _{\RR} := V \otimes _{\QQ} \RR = K _{\RR} ^{n}$
and for any $x\in V_{\RR }$ the matrix 
$x x^{\dagger } $ lies in $\Sigma  $.
The following lemma is easily checked :
\begin{lemma}
Any $F\in \Sigma  $ defines a 
quadratic form on $V_{\RR }$ by :
 $$F[x] := \langle F, x x^{\dagger } \rangle \in \RR  \mbox{ for all } x\in V_{\RR} .$$
This quadratic form is positive definite if and only if $F\in {\mathcal P}$.
\end{lemma}
As a consequence, with a slight abuse of language, we will sometimes refer to elements of $\Sigma$ as \emph{forms}.

\section{Minimal vectors}\label{minvec} 
Let $A=M_n(K)$ for some division algebra $K$. As before we fix some maximal order ${\mathcal O}$ in $K$ and choose some 
right ${\mathcal O}$-lattice $L$ in the simple left $A$-module $V = K^n$. 
Then $\Lambda := \End _{\mathcal O}(L)$ is a maximal order in $A$ with 
unit group $\Lambda ^{\times} := \GL(L) = \left\lbrace  a\in A \mid a L = L \right\rbrace  $. 

Following \cite{Ash84}, we will define the $L$-minimum of a form $F \in \mathcal P$ with respect to a \emph{weight}. 
\begin{definition} A weight $\varphi$ on $L$ is a  $\GL(L)$-invariant map from the 
projective space $\mathbf{P} (K^n)$ to the positive reals,
 such that $\max_{x \in \mathbf{P} (K^n)}\varphi(x)=1$.
\end{definition}

A natural choice for the weight is $\varphi_0(x)=1$ for all 
$x \in K ^n -\left\lbrace 0 \right\rbrace$. 
However, another rather standard choice for $\varphi$ is possible,
 which allows for definitions having a natural geometric interpretation
 and somehow simplify the computations, at least  in the case of imaginary quadratic fields
 or definite quaternion algebras (see Section \ref{mini}).
 Roughly speaking, this alternative weight is given by the inverse of the
 $\gcd$ of the coefficients of a vector in $K^n$ with respect to a given pseudo-basis
 of the lattice $L$. To be more precise, we need the following definition
\begin{definition}
Let $L=e_1 \mathfrak{c}_1  \oplus ... \oplus e_n \mathfrak{c}_n$. To any $\ell = \sum _{i=1}^n e_i \ell _i \in L -\left\lbrace 0\right\rbrace$ we associate the 
integral left $\OO$-ideal 
$${\mathfrak a}_{\ell } := \sum _{i=1}^n  \mathfrak{c}_i^{-1} \ell _i$$
as well as its norm
$$\N( {\mathfrak a}_{\ell } ) := |\OO / {\mathfrak a}_{\ell } | =\N_{Z(K)/\QQ}\left(\nr\left(  {\mathfrak a}_{\ell }\right)^d\right) .$$
\end{definition}

\begin{lemma}\label{weight} 
\begin{enumerate}
\item[(a)] $\N( {\mathfrak a}_{\ell } ) \geq 1$ for all $\ell \in L-\left\lbrace 0\right\rbrace$.
\item[(b)] For any $\lambda \in K^{\times}$ and $\ell = \sum _{i=1}^n e_i \ell _i \in L-\left\lbrace 0\right\rbrace$, one has ${\mathfrak a}_{\ell \lambda } = \a _{\ell } \lambda$.
\item [(c)] If $g\in \GL(L)$ and  $\ell = \sum _{i=1}^n e_i \ell _i \in L-\left\lbrace 0\right\rbrace$, 
then ${\mathfrak a}_{g\ell } = \a _{\ell }$. 
\end{enumerate}
\end{lemma}

\bew
(a) is clear, because all $\mathfrak{c}_i^{-1} \ell _i$ are integral left $\OO $-ideals, and (b) is straightforward.
\\
 To see (c) 
 write $g e_i = \sum _{j=1}^{n} e_j g_{ji }$.
 Since $gL \subseteq L$ we get  
$g_{ji}  \in {\mathfrak c}_j {\mathfrak c}_i^{-1} $. 
Then 
$ g \ell =  \sum _{j=1}^n e_j (\sum _{i=1}^{n} g_{ji} \ell _i )  $ and
$$\a _{g\ell } = \sum _{j=1}^n {\mathfrak c}_j^{-1} \sum _{i=1}^n g_{ij} \ell _i 
\subseteq \sum _{j,i} {\mathfrak c}_j^{-1} {\mathfrak c}_j {\mathfrak c}_i^{-1} \ell _i \subseteq \a _{\ell }  .$$
One obtains equality by applying $g^{-1} \in \GL(L) $. 
\eb

Now for any $x \in K^n$, we can find $\lambda \in K-\left\lbrace 0\right\rbrace$ 
such that $x \lambda \in L$. 
It follows from the previous lemma that the class of
 $\nr(\a_{x \lambda})$ in $\Cl _K(R)$ does not depend on the choice of
 an element $\lambda$ with this property. 
Consequently, if we define the norm of a class in $\Cl _K(R)$ 
as the smallest possible norm of an integral ideal in that class, 
we can associate to $x$ a well-defined quantity $N_x$ by the formula
$$N_x=\N( \left[ \nr\left(\a_{x \lambda}\right) \right])=\min\limits_{\substack{I \subset \OO  \\ \left[ \nr\left(I\right)\right] = \left[ \nr\left(\a_{x \lambda}\right) \right]}} \N_{Z(K)/\QQ}\left(\nr\left(I\right)^d\right),$$
where as before $\lambda$ is any element in $K-\left\lbrace 0\right\rbrace$ such that $x \lambda \in L$. 
This in turn can be used to define a weight $\varphi_1$ on $K^n$
setting 
\begin{equation}\label{fi1} 
\varphi_1(x)=N_x^{-2/[K:\QQ]}
\end{equation}
(that this is indeed a weight follows immediately from Lemma \ref{weight}).

\begin{remark}
As explained in \cite{Ash84}, the space of weights is isomorphic to $\RR^{h_K-1}$, where $h_K$ stands for the class number of $K$. In particular, the trivial weight $\varphi_0$ is the only possible choice if $h_K=1$ (and $\varphi_1=\varphi_0$ in that case).
\end{remark}

Having fixed a weight $\varphi$ on $L$, we can define the minimum of a form an its set of minimal vectors as follows :
\begin{definition}\label{defallg}
The $L$-minimum of $F \in {\mathcal P}$ with respect to the weight $\varphi$ is 
\begin{equation*}\min\nolimits _L (F) := \min_{\ell \in L-\left\lbrace 0\right\rbrace } \varphi(\ell) F[\ell ]. \end{equation*}
The set of minimal vectors of $F$ in $L$ is defined as
\begin{equation*} S_L(F) := \left\{ \ell \in L - \left\lbrace 0\right\rbrace  \mid \varphi(\ell) F[\ell] = \min\nolimits _{L}(F) \right\}.  \end{equation*}
\end{definition}

\begin{remark}
The set $S_L(F)$ is finite. Indeed, let 
$m:= \min  \left\lbrace \varphi(\ell) \mid \ell \in L\setminus \left\lbrace 0 \right\rbrace \right\rbrace   $. Then $m>0$ as $\varphi $ takes only finitely many positive real values, so
 $S_L(F) \subset \left\lbrace  \ell \in L \mid F[\ell ] \leq m^{-1} \min\nolimits _L (F)\right\rbrace  $ which is a finite set and can be computed as the set of vectors of small length
in a $\ZZ $-lattice.
\end{remark}

\section{Minimal classes}\label{minclass} 
We keep the general assumptions of the previous section:  $K$ is a division algebra,
 ${\mathcal O}$ a maximal order in $K$ and  $L$ a right ${\mathcal O}$-lattice in $K^n$,
 on which a weight $\varphi$ is fixed.
\begin{definition}
Two elements $F_1$ and $F_2 \in {\mathcal P}$ are called
\emph{minimally equivalent with respect to $L$ and $\varphi$}, if $S_L(F_1) = S_L(F_2) $.
We denote by
$$\Cl _L (F_1) := \left\lbrace  F \in {\mathcal P} \mid S_L(F) = S_L(F_1) \right\rbrace  $$
the \emph{minimal class} of $F_1$.
If $C=\Cl _L (F_1)$ is a minimal class then we define
$S_L(C) = S_L(F_1)$ the associated set of minimal vectors.
A minimal class $C$ is called \emph{well rounded}, if $S_L(C)$ contains a
$K$-basis of $V$.
The form $F\in {\mathcal P}$ is called \emph{perfect with respect to $L$},
if $\Cl _L(F) = \left\lbrace  a F \mid a \in \RR , a > 0 \right\rbrace  $.
\end{definition}

\begin{remark} Note that minimal classes and all subsequent definitions in this section 
actually \emph{depend on} 
 the weight $\varphi$, although we do not indicate it systematically
 in our notations. No inconstancy can arise from this, 
since we work with fixed weight $\varphi$ (and fixed lattice $L$).
\end{remark}

The group $\GL_n(K)$, and hence also its subgroup $\GL(L)$, acts on $\Sigma $ by 
$(F,g) \mapsto g^{\dagger } F g  $ (where we embed $A$ into $A_{\RR} $
to define the multiplication). 
Two forms in $\Sigma $ are called \emph{$L$-isometric}, 
if they are in the same $\GL(L)$-orbit.
For $F\in {\mathcal P}$ we denote by 
$$\Aut _L (F) :=  \left\lbrace  g\in \GL(L) \mid g^{\dagger} F g = F \right\rbrace  $$
the \emph{automorphism group} of $F$. 
Then
$\Aut_L(F) $ is always a finite subgroup of $\GL (L)$. 
The group $\GL(L)$ acts on the set of minimal classes.
Two minimal classes are called \emph{equivalent}, if they are in the
same orbit under this action. The stabiliser of a minimal class is
called the \emph{automorphism group} of the class,
$$\Aut _L (C) = \left\lbrace  g\in \GL(L) \mid g S_L(C) = S_L(C) \right\rbrace  .  $$

\begin{lemma} (see \cite[Proposition 2.9]{Batut})\label{aut} 
Let $C$ be a well rounded minimal class.
Then the \emph{canonical form}
$T_C:= \sum _{x\in S_L(C)} x x^{\dagger } \in {\mathcal P} $ is positive definite and
$\Aut_L (C) = \Aut _L(T_C^{-1})$.
Two well rounded
minimal classes $C$ and $C'$ are equivalent, if and only
if $T_C^{-1}$ and $T_{C'}^{-1}$ are $L$-isometric.
\end{lemma}

\bew
The proof is similar to the one in \cite{Batut}.
The well roundedness of $C$ implies that the rank of $T_C$ is maximal:
The mapping  $(,): V\times V \to K_{\RR }, 
(x,y) := x^{\dagger } y $ is Hermitian and 
non-degenerate. 
Let $\left\lbrace  x_1,\ldots , x_n \right\rbrace  \subset S_L(C)$  be a  $K$-basis of $V$,
then for any $v\in V$
$$\sum _{i=1}^n x_i x_i^{\dagger } v = \sum _{i = 1}^n  x_i ( x_i, v) = 0 
\mbox{ if and only if } v\in V^{\perp } = \left\lbrace  0 \right\rbrace  $$
so the kernel of the positive semidefinite matrix $\sum _{i=1}^n x_i x_i^{\dagger} $
is $\left\lbrace  0 \right\rbrace  $, therefore $T_C $ is invertible and hence in ${\mathcal P}$.
Clearly $\Aut_L (C)
 \subseteq \Aut _L(T_C^{-1})$.
To see the converse
put  $s:= |S_L(C)|$
and let $S \in M_{n,s}(K)$ be a matrix whose columns are the elements
of $S_L(C)$, in particular   $T_C = S S^{\dagger }$.
Take some $g\in \Aut(T_C^{-1}) = \left\lbrace  g\in \GL(L) \mid g T_C g^{\dagger } = T_C \right\rbrace $
and put $S' := gS $.
Then $S'(S')^{\dagger } = T_C = S S^{\dagger }  $
and for any $F\in {\mathcal P} $
$$(\star ) \ \ \sum _{y\in \cols( S') } F [y ]  = \trace ((S')^{\dagger } F S') = 
\langle S'(S')^{\dagger } ,F  \rangle = \langle S S^{\dagger } , F \rangle  =
\sum _{x\in S_L(C) } F [ x ] .$$
If $x$ is some column of $S$ and $y:=gx$, then $\varphi(y)=\varphi(x)$, because of the $\GL(L)$-invariance of $\varphi$. Moreover $\varphi(y)F[y] \geq \varphi(x)F[x] =\min_{\ell \in L - \left\lbrace 0\right\rbrace }\varphi(\ell)F[\ell]$, whence $F[y] \geq F[x]$, with equality if and only if $y\in S_L(C)$. So we can only have equality in $(\star )$
 if $S_L(C) = \left\lbrace  \cols (S' )\right\rbrace  $ and hence $g\in \Aut_L(C) $.
\eb

\section{Maximal finite subgroups of $\GL(L)$}
In this section we will use variants of the Voronoi algorithm
to compute a set of representatives of the conjugacy classes
of maximal finite subgroups of $\GL(L)$.
The known methods (see e.g. \cite{PP76}) start 
with  the list of all conjugacy classes of finite subgroups
of $\GL_n(K)$. For each group $G$ they compute the invariant lattices 
to find the $\GL(L)$-conjugacy classes of subgroups in the class of $G$.
In particular for reducible groups $G$ this set of invariant lattices is infinite
and one needs to use the action of $N_{\GL_n(K)}(G)$. 
Also it seems to be difficult to restrict to one isomorphism class of $\OO $-lattices $L$. 

Here we will start with some lattice $L$ and use the tessellation of 
the cone of positive definite hermitian forms into $L$-minimal classes 
to obtain a list of subgroups of $\GL(L)$ that contains representatives of 
all conjugacy classes of maximal finite subgroups of $\GL(L)$. 
To check maximal finiteness and also conjugacy of the groups in
the list, we use a relative version of Voronoi's theory. 

\begin{definition}
Let $G\leq \GL (L)$ be some finite subgroup. 
Let ${\mathcal F}(G) := \left\lbrace  F\in \Sigma  \mid  g^{\dagger } F g = F  \mbox{ for all } g\in G \right\rbrace  $
 denote the space of $G$-invariant Hermitian forms.
It contains the cone
${\mathcal F}_{>0} (G) := {\mathcal F}(G) \cap  {\mathcal P}  $ 
of positive definite 
$G$-invariant forms. 
For $F\in {\mathcal F}_{>0} (G) $ the 
\emph{$G$-minimal class} of $F$ is $\Cl_L(F) \cap {\mathcal F}(G)$.
A form $F\in {\mathcal F}_{>0} (G) $ is called
 \emph{$G$-perfect with respect to $L$}, 
if 
 $\Cl_L(F) \cap {\mathcal F}(G) = \left\lbrace  aF \mid a\in \RR _{>0} \right\rbrace  $.
\end{definition}

\begin{lemma}\label{Ginvmin}
Let 
$$\pi _G: \Sigma  \to {\mathcal F}(G) , 
F \mapsto \frac{1}{|G|} \sum _{g\in G} g^{\dagger } F g $$ 
be the usual averaging operator and $C$ be a $G$-invariant minimal class. 
Then $$C\cap {\mathcal F}(G)  = \pi _G(C).$$
\end{lemma}

\bew
Since $\pi _G(F) = F $ for all $G$-invariant forms, it is clear that 
$C\cap {\mathcal F}(G)  \subseteq \pi _G(C) $.
So let $F\in C$. Then $S_L(F) = S_L(C) $. 
Since $S_L(C)$ is $G$-invariant, $S_L(C) = S_L(g^{\dagger }F g) $ for any $g\in G$.
As $\pi _G(F)$ is a sum of positive definite forms, also 
$S_L(\pi _G(F)) = S_L(C)$ and so $\pi _G(F) \in C$.
\eb

As in the classical case, Voronoi's algorithm, as described e.g. in \cite{Opgenorth} can be adapted to the case of $G$-invariant forms to 
compute the $G$-perfect forms and the cellular decomposition of 
${\mathcal F}_{>0}(G)$ into $G$-minimal classes up to the 
action of the normaliser (see for instance \cite[Theorem 2.4]{Batut} for details on this procedure in the classical case). 

\begin{proposition}
Let $G\leq \GL(L)$ be finite. Then 
there exists at least one $G$-perfect form with respect to $L$.
\end{proposition}

\bew
We will show that $L- \left\lbrace 0 \right\rbrace $ is discrete and admissible 
in the sense of \cite[Definition 1.4]{Opgenorth}.
Then by \cite[Proposition 1.8]{Opgenorth} 
there exists a $G$-perfect form. Moreover, \cite[Theorem 1.9]{Opgenorth} tells us that the Voronoi domains of the $G$-perfect forms
form an exact tessellation of ${\mathcal F}(G^{\dagger}) _{>0 }$.
\\
Clearly $L- \left\lbrace  0 \right\rbrace $ is discrete in $V_{\RR} := V\otimes _{\QQ} \RR $. 
For the admissibility we need to show that for any $F\in \partial {\mathcal P}$,
the boundary of ${\mathcal P}$, and any $\epsilon >0$ there is some 
$\ell \in L - \left\lbrace  0 \right\rbrace $, such that $\varphi(\ell ) F[\ell ] < \epsilon $.
If $F\in \partial {\mathcal P}$, it is positive semidefinite,
so 
$$\left\lbrace  x\in V_{\RR } \mid  F[x] = 0 \right\rbrace  
 = \left\lbrace  x\in V_{\RR } \mid Fx = 0 \right\rbrace  \leq V_{\RR} $$
is a subspace. 
In particular the volume of the  convex set 
$$ {\mathcal K} _{\epsilon } := \left\lbrace   x\in V_{\RR } \mid 
 F[x]  < \epsilon \right\rbrace   = - {\mathcal K} _{\epsilon } $$
is infinite and by Minkowski's convex body theorem
${\mathcal K}_{\epsilon } $ contains some $0\neq \ell \in L $. 
Then $F[\ell ]< \epsilon $ and hence also 
$\varphi(\ell ) F[\ell ] < \epsilon $ since $\varphi(\ell ) \leq 1$.
\eb

\begin{lemma}\label{wr}
Let $G\leq \GL(L)$ be finite.
Then any $G$-perfect form $F$ with respect to $L$ is well rounded. 
\end{lemma}

\bew
The proof is similar to the classical case. 
Assume that  $\langle S_L(F) \rangle _K \neq V $. 
Then there is some linear form $H\in V^* = K^n$ so that $Hx= 0 $ 
for all $x\in S_L(F)$. 
Let 
$$F_0 := \frac{1}{|G|}  \sum _{g\in G} g^{\dagger } H^{\dagger } H g .$$
Since $S_L(F)$ is $G$-invariant, $x^{\dagger }F_0x = 0$ for all $x\in S_L(F)$,
so $S_L(F+\epsilon F_0) \supseteq S_L(F) $ for all $\epsilon > 0$
with equality, if $\epsilon $ is small enough. 
So $F+\epsilon F_0 \in \Cl _L (F) \cap {\mathcal F}_{> 0} (G) $ 
contradicting the assumption that $F$ is $G$-perfect with respect to $L$.
\eb

\begin{theorem} \label{minmax}
Let $G \leq \GL(L)$ be some maximal finite subgroup of $\GL(L)$.
Then 
$G=\Aut _L (C)$ 
for some  well rounded minimal class $C$ with respect to $L$,
 such that 
 $ C \cap {\mathcal F} (G)$ spans a subspace of ${\mathcal F}(G)$ of dimension 1.
\end{theorem}

\bew
The group $G$ always fixes some $G$-perfect form $F$ with respect to $L$. 
Let $C:= \Cl _L (F)$. Then $S_L(C) = S_L(F)$ is $G$-invariant, 
so $G\leq \Aut _L(C)$. 
By Lemma \ref{wr} $C$ is well rounded, so $\Aut _L(C)$
is finite and the maximality of $G$ implies that $G=\Aut _L(C)$.
\eb

With Theorem \ref{minmax}  we obtain a finite list of 
finite subgroups of $\GL(L)$ that contains a system of representatives 
of conjugacy classes of maximal finite subgroups. 
We need to be able to 
test maximal finiteness and conjugacy in $\GL(L)$ of such groups
$\Aut_L(C) $.

\begin{proposition} \label{obermin}
Let $G\leq \GL(L)$ be some finite subgroup. 
Then the maximal finite subgroups $H$ of $\GL(L)$ that contain $G$
are of the form $H = \Aut_L(C_G)$ for some $G$-minimal class $C_G$.
\end{proposition} 

\bew
Let $H$ be some maximal finite subgroup of $\GL(L)$ that contains $G$.
By Theorem \ref{minmax} 
 there is some $G$-invariant $L$-minimal class $C$ such that 
$H=\Aut_L(C)$.
By Lemma \ref{Ginvmin} $S_L(C)= S_L(C_G)$ for the $G$-minimal class 
$C_G =\pi_G(C)$ 
and $H = \Aut_L(C_G) $.
\eb

\begin{remark}
To test whether two maximal finite subgroups $G_1$, $G_2$
 of $\GL(L)$ are conjugate, 
one computes a system of representatives $R_i$ of the 
$N_{\GL(L)}(G_i)$-orbits of $G_i$-perfect forms and then 
checks whether a given form in $R_1$ is $L$-isometric to
some form in $R_2$. 
Since $G_i = \Aut_L(F_i) $ for all $F_i\in R_i$, any isometry 
yields a conjugating element. 
\end{remark}

\section{Imaginary quadratic fields and definite quaternion algebras}\label{mini}
In this section we will assume that 
$K$ is either the field of rational numbers, an imaginary quadratic number field or a
definite quaternion algebra over the rationals.
These are exactly the cases, where $K_{\RR }$ is a division algebra and 
$\Sym (K_{\RR} ) = \RR $. We thus have in those cases (and in those cases only) the noteworthy property that 
\begin{equation}\label{center} 
\forall \lambda \in K_{\RR }, \forall x \in V_{\RR } \quad  F[x \lambda]= \lambda \lambda^{\dagger } F[x].
\end{equation}
As a consequence, it is more natural and more efficient to compute minima with respect to the weight $\varphi_1$ defined in the previous subsection, because of the following lemma
 \begin{lemma}\label{cusp} 
 For any $F \in {\mathcal P}$ one has
\begin{equation*}\min\nolimits _L (F) := \min_{\ell \in L-\left\lbrace 0\right\rbrace } \frac{F[\ell]}{\N( {\mathfrak a}_{\ell } )^{2/[K:\QQ]}} \end{equation*}
where the minimum on the left hand side is computed with respect with the weight $\varphi_1(\ell)=N_{\ell}^{-2/[K:\QQ]}$.
 \end{lemma}
\bew  The inequality $\min\nolimits _L (F) \geq \min_{\ell \in L-\left\lbrace 0\right\rbrace } \frac{F[\ell ]}{\N( {\mathfrak a}_{\ell } )^{2/[K:\QQ]}}$ is clear, since $N_{\ell} \leq \N( {\mathfrak a}_{\ell })$ for every $\ell \in L-\left\lbrace 0\right\rbrace $. In the opposite direction, every $\ell \in L-\left\lbrace 0\right\rbrace $, there exists $\lambda \in K - \left\lbrace 0\right\rbrace $ such that ${\mathfrak a}_{\ell\lambda}={\mathfrak a}_{\ell}\lambda \subset \OO$ and $\N({\mathfrak a}_{\ell}\lambda)=\N([{\mathfrak a}_{\ell }])=N_{\ell}$ (in particular, $\ell \lambda \in L$) . Using (\ref{center}), we see that $ \frac{F[\ell]}{\N( {\mathfrak a}_{\ell } )^{2/[K:\QQ]}}= \frac{F[\ell \lambda]}{\N( {\mathfrak a}_{\ell \lambda})^{2/[K:\QQ]}}=\varphi_1(\ell \lambda) F[\ell \lambda] \geq \min_L(F)$, whence the conclusion taking the minimum of the left-hand side over $ L-\left\lbrace 0\right\rbrace $.
\eb

\begin{remark} The reformulation given in Lemma \ref{cusp} of the minimum of a form with respect to $\varphi_1$ has two noteworthy applications
\begin{enumerate}
\item It can be interpreted in terms of  \emph{minimal distance to cusps} as explained in \cite{Men79} (see also \cite[chapter 7]{EGM98}). 
\item One can easily deduce from this that the Voronoi complex will depend only on the Steinitz class of L modulo $n$th powers (see \cite[Theorem 3.8]{Braun}).
\end{enumerate}

\end{remark}

\section{Examples}\label{exples} 

We will use the method from the previous section to 
compute the conjugacy classes of maximal finite subgroups of 
$\GL(L)$ for imaginary quadratic number fields $K$. 
This is an invariant of the isomorphism class of $\GL(L)$ and 
will show that for small examples 
these groups  are not isomorphic.

\noindent
{\bf Example 1} \\
Let $K:=\QQ [\sqrt{-15}]$, $\OO = O_K = \ZZ [\frac{1+\sqrt{-15}}{2} ]$, $n=2$. Then $\Cl(K) = \left\lbrace  [O_K] , [\wp _2] \right\rbrace $ 
where $\wp _2$ is some prime ideal dividing $2$, 
so there are two isomorphism classes of $O_K$-lattices in $K^2$: one corresponding to the lattice
$L_0$ with Steinitz-invariant $[O_k]$ and the other one to the lattice $L_1$ with Steinitz-invariant $[\wp _2]$. 
For both lattices the perfect forms are given in \cite{Braun}. 
\\
For both lattices $L$, Table 1 lists the $\GL (L) $-orbits of well rounded
minimal classes $C$ 
according to their perfection corank  together with their stabilizers 
$G=\Aut_L(C)$. The two classes of perfection corank 0 contain the perfect forms.
The third column gives the dimension of $\pi _G(C) $.
If this dimension is one, then $\pi_G(C) \subset  \langle F \rangle $ 
for some $G$-perfect form $F$,  the next column gives the 
automorphism group $\Aut_L(F)$ and the last one indicates whether $G$ is maximal finite. 

\begin{table}[pbth]
\caption{Well rounded minimal classes for $K=\QQ [\sqrt{-15}]$}
\centering
\begin{tabular}{|c|c|c|c|c|}
\hline
\multicolumn{5}{|c|}{$L=L_0$} \\
\hline
$C$ & $G=\Aut_L(C) $ & $\dim  (\pi_G(C) )$ & $\Aut_L(F)$  & maximal \\
\hline
\multicolumn{5}{|c|}{perf. corank = 0} \\
\hline
$P_1  $ & $C_6 $ & 1 & $C_6$& no \\
$P_2  $ & $C_4 $ & 1 & $C_4$& no \\
\hline
\multicolumn{5}{|c|}{perf. corank = 1} \\
\hline
$C_1$ & $D_{12}$ &  1 & $ D_{12} $ & yes \\
$C_2$ & $D_{12}$ &  1 & $ D_{12} $ & yes \\
$C_3$ & $C_{2}$ &  2 & & no  \\
$C_4$ & $C_{2}$ &  2 & & no  \\
\hline
\multicolumn{5}{|c|}{perf. corank = 2} \\
\hline
$D_1$ & $D_{8}$ &  1 & $ D_{8} $ & yes \\
$D_2$ & $D_{8}$ &  1 & $ D_{8} $ & yes \\
$D_3$ & $V_{4}$ &  1 & $ V_{4} $ & yes \\
$D_4$ & $V_{4}$ &  1 & $ V_{4} $ & yes \\
\hline
\multicolumn{5}{|c|}{$L=L_1$} \\
\hline
$C$ & $G=\Aut_L(C) $ & $\dim  (\pi_G(C) )$ & $\Aut_L(F)$ & maximal \\
\hline
\multicolumn{5}{|c|}{perf. corank = 0} \\
\hline
$P  $ & $C_3:C_4 $ & 1 & $C_3:C_4$ & yes \\
\hline
\multicolumn{5}{|c|}{perf. corank = 1 }\\
\hline
$C_1$ & $D_{8}$ &  1 & $ D_{8} $  & yes \\
$C_2$ & $D_{8}$ &  1 & $ D_{8} $  & yes \\
$C_3$ & $D_{12}$ &  1 & $ D_{12} $  & yes \\
\hline
\multicolumn{5}{|c|}{ perf. corank = 2} \\
\hline
$D$ & $V_{4}$ &  1 & $ V_{4} $  & yes \\
\hline
\end{tabular}
\end{table}

\noindent \underline{$L=L_0$}: \\
The two groups $G=D_8$ and $G=D_{12}$ are absolutely irreducible
maximal finite subgroups of $\GL_2(K)$. Since $\dim ({\mathcal F} (G)) = 1$ for both groups
and $C_i$ and $D_i$ are inequivalent $(i=1,2)$
one gets 2 conjugacy classes of maximal finite subgroups of both 
isomorphism types $D_8$ and $D_{12}$. 
To prove that $G:=\Aut_L(D_3)$ is maximal finite, we
compute the well rounded $G$-minimal classes, using Voronoi's algorithm
and starting with the $G$-perfect form $F\in \pi_G(D_3)$.
$S_L(F) = \left\lbrace  \pm v_1,\pm v_2 \right\rbrace  $ with $\a _{v_1} = O_K$, $\a_{v_2} = \wp _2$.
Therefore both minimal vectors are $G$-eigenvectors and the 
$G$-Voronoi domain has 2 faces, both of which are dead ends (see \cite[Definition 13.1.7]{Mart}). 
So $F$ is the unique $G$-perfect form and there are no other 
well rounded $G$-minimal classes. 
The situation is the same for $\Aut_L(D_4)$. The two $G$-perfect 
forms in $D_3$ and $D_4$
(rescaled to have minimum 1) are Galois conjugate 
but not $L$-isometric, with shows that 
$\Aut_L(D_3)$ and $\Aut_L(D_4)$ are  not conjugate in $\GL(L)$. 

The proof that $G:=\Aut_L(P_i)$ is
not maximal finite is similar for both cases $i=1,2$. 
The space of invariant forms has dimension 2, there
are two $G$-orbits on $S_L(P_i)$, so the $G$-Voronoi domain of
$P_i$ has two faces, corresponding to 1-dimensional $G$-minimal classes 
with automorphism group $D_{12}$ (for $P_1$) resp. $D_8$ (for $P_2$). 
One checks for $i=1,2$ that  $\Aut_L(P_i) $ is properly contained in these groups.

\noindent \underline{$L=L_1$:} 
As in the free case the uniform groups $\Aut_L(P)$ and $\Aut_L(C_i)$, 
$i=1,2,3$ are maximal finite and represent distinct conjugacy classes.
For the group $G = \Aut_L(D) \cong V_4$ we again have a 
unique $G$-perfect form $F$ and the two $L$-minimal vectors of $F$ 
are eigenvectors for $G$. So both faces of the $G$-Voronoi domain of 
$F$ are dead ends and $G=\Aut_L(F)$ is maximal finite. 

As $\GL_2(O_K)$ and $\GL(L_1)$ have different conjugacy classes of 
maximal finite subgroups one finds the following corollary.

\begin{folg} 
$\GL_2(O_K) = \GL(L_0)$ and $\GL(L_1)$ are not isomorphic. 
\end{folg} 

\begin{table}[pbth]
\caption{Number of conjugacy classes of maximal finite subgroups}
\centering
\begin{tabular}{|c|c|c|c|c|c|c|}
\hline 
 & $  D_8 $  & $  D_{12} $  & $  V_4 $  & $  \SL_2(3) $  & $  Q_8 $  & $  C_3:C_4 $  \\ 
\hline
 $ K = \QQ[\sqrt{-15}] $  & & & & & & \\
 $ St(L) = [O_K ] $  & 2& 2& 2& - & - & - \\
 $ St(L) = [\wp _2 ] $  & 2& 1& 1& - & - & 1 \\
\hline
 $ K = \QQ[\sqrt{-5}] $  & & & & & & \\
 $ St(L) = [O_K ] $  & 3& 2& 1& - & 1 & - \\
 $ St(L) = [\wp_2 ] $  & 1& 2& 1& 1 & - & - \\
\hline
 $ K = \QQ[\sqrt{-6}] $  & & & & & & \\
 $ St(L) = [O_K ] $  & 3& 2& 1& 1 & - & - \\
 $ St(L) = [\wp_2 ] $  & 1& 1& 2& - & 1 & 1 \\
\hline
 $ K = \QQ[\sqrt{-10}] $  & & & & & & \\
 $ St(L) = [O_K ] $  & 3& 2& 1& - & 1 & - \\
 $ St(L) = [\wp_2 ] $  & 1& -& 3& 1 & - & 2 \\
\hline
 $ K = \QQ[\sqrt{-21}]  $ & & & & & & \\
 $ St(L) = [O_K ] $  & 6& 4& 2& - & - & 2 \\
 $ St(L) = [\wp_2 ] $  & 2& -& 6& - & - & - \\
 $ St(L) = [\wp_3 ] $  & -& 2& 6& 2 & - & - \\
 $ St(L) = [\wp_5 ] $  & -& -& 8& - & 2 & - \\
\hline 
\end{tabular}
\end{table}

\noindent {\bf Example 2}  \\
Table 2 lists the results of 
similar computations which we did for 
 certain small imaginary quadratic fields. In particular we find

\begin{corollary}
Let $K$ be one of the six fields in Table 2.
Then non-conjugate maximal orders in $M_2(K)$ have
non-isomorphic unit groups.
\end{corollary}


\end{document}